\documentclass[reqno,12pt]{amsart}
\usepackage{enumitem}
\usepackage{stmaryrd}
\usepackage{amsmath}
\usepackage{amssymb}
\usepackage{amsthm}
\usepackage[english]{babel}

\newtheorem{theorem}{Theorem}

\newtheorem{lemma}{Lemma}

\begin{document}
\title[Transport distance estimates]{Estimates of transport distance \\ \noindent in the central limit theorem}

\author{Andrei~Yu.~Zaitsev}

\address{St. Petersburg Branch of the Steklov Mathematical Institute
\\
Fontanka 27 \\ St. Petersburg 191023, Russia\\
and St. Petersburg State University, \\ Universitetskaya Embankment 7/9 \\ St. Petersburg,
199034 Russia}
\email{zaitsev@pdmi.ras.ru}

\keywords {inequalities, sums of independent random vectors, transport distance estimation, central limit theorem}

\maketitle

\section{Introduction}

Sums of independent terms first appeared in probability theory when considering binomial distributions within the Bernoulli scheme. The law of large numbers and the Moivre-Laplace central limit theorem were derived.
It was observed that binomial distributions are not only well approximated by the normal law, but also that the tail decay of binomial distributions is similar to that of normal distributions. As a natural extension of the class of binomial distributions, we can consider the class of distributions of sums of independent (generally, non-identically distributed) random variables bounded in absolute value by the same constant.
A~significant number of works are devoted to estimating the tails of such distributions, see, for example, \cite{benn, b02, b04, H,N,p14},
In this paper we will discuss not only the tail decay, but also the estimation of the approximation accuracy in the one-dimensional central limit theorem.

Let $X,X_1,\ldots,X_n$ be $d$-dimensional independent random vectors bounded with probability one. For simplicity, we will assume that they have zero mean values:
\begin{equation}\label{tau}
\mathbf{P}\{\|X_{j}\|\le\tau\}=1,\quad\mathbf{E}\,X_{j}=0,\quad j=1,\ldots, n.
\end{equation}
We will be interested in the behavior of the distribution of the sum $S=X_{1}+\cdots+X_{n}$ depending on the limiting value $\tau>0$.

From the non-uniform Bikelis estimate \cite{Bik} in the one-dimensional central limit theorem it follows that in the one-dimensional case
\begin{equation}\label{bik}
  W_1(F,\Phi)\le c\tau.
\end{equation} with an absolute constant $c$,
where $W_1$ is the Kantorovich--Rubinstein--Wasserstein transport distance (see review articles \cite{BK,B22}), $F=\mathcal{L}(S)$ is the distribution of the sum $S$, and $\Phi=\Phi_F$ is the corresponding normal distribution, with the same zero mean and the same variance as those of the distribution $F$.
When proving  inequality \eqref{bik}, it should be taken into account that, according to \cite{V},
\begin{equation}\label{val}
  W_1(F,\Phi)=\int|F(x)-\Phi(x)|\,dx,
\end{equation}where
$F(\,\cdot\,)$ and $\Phi(\,\cdot\,)$ are the corresponding distribution functions. In addition, $\mathbf{E}\,|X_{j}|^{3}\le\tau\,\mathbf{E}\,X_j^{2}$.

The main
result of the paper is significantly stronger and more precise. It is claimed that
\begin{equation}\label{thm}
 W(F,\Phi) =\inf_{\pi}\int\exp\{|x-y|/c\tau\}\,d\pi(x,y)\le c,
\end{equation}
where the infimum is taken over all two-dimensional probability distributions $\pi$ with marginal distributions $F$ and $\Phi$. The result is also generalized to distributions with sufficiently slowly growing cumulants from the class $\mathcal{A}_{1}(\tau )$ introduced in the author's paper~\cite{hh29}. In special cases, we obtain some results of Rio~\cite{Rio}. The possibility of generalizing the result to the multidimensional case is discussed.

Following Rio~\cite{Rio}, we define the Wasserstein distance associated with the Orlicz function~$\psi$:
\begin{equation}\label{th}
 W_\psi(G,H) =\inf\Big\{a>0:\inf_{\pi}\int\psi(|x-y|/a)\,d\pi(x,y)\le 1\Big\},
\end{equation}
where the second infimum is taken over all two-dimensional probability distributions $\pi$ with marginal distributions $G$ and $H$.

Inequality \eqref{thm} can be rewritten as
\begin{equation}\label{thm3}
 W_\psi(F,\Phi)\le c\tau,
\end{equation}
with the Orlicz function $\psi(x)=\exp\{|x|\}-1$. Inequality \eqref{bik} can also be written in the form~\eqref{thm3}, but for the Orlicz function $\psi(x)=|x|$.
Inequality~\eqref{thm3} is also valid for the Orlicz function $\psi(x)=|x|^{p}$, $p\ge1$. In this case, the statement is easily deduced from~\eqref{thm3} and turns into the estimate
\begin{equation}\label{bikp}
  W_p(F,\Phi)\le c(p)\,\tau,
\end{equation}
where $W_ P(\,\cdot\,,\,\cdot\,)$ is the standard Wasserstein $p$-distance. We took into account that $|x|^{p}\le c(p)\,\exp\{|x|\}$.

The class of distributions of sums $S=X_{1}+\cdots+X_{n}$ satisfying  conditions \eqref{tau} can be considered as a natural generalization of the class of binomial distributions, which historically turned out to be the first distributions of sums of independent terms studied in probability theory. Bernstein~\cite{Bern} found less restrictive conditions (see  definition~\eqref{bern} for $d=1$), under which the tails of the distributions of sums admit estimates similar to those for the tails of binomial distributions. Under the conditions of Bernstein's inequality, the distributions of the terms have finite exponential moments, that is, the Cram\'er conditions are satisfied, under which the theorems on large deviations of the distributions of sums of independent terms are valid. As is well known, the coefficients of the so-called Cram\'er--Petrov series arising in the formulations are determined from the cumulants of the distributions of the sums, see \cite[Lemma 1.4]{91}. This motivated Statulevi\v {c}ius~\cite{hh27} to a further expansion of the class of distributions for which the results on large deviations are valid. He introduced classes of distributions that are no longer necessarily representable as distributions of sums of a large number of independent terms, but whose cumulants behave similarly to the cumulants of such sums (see~\eqref{stat}). In this paper, we prove  inequality~\eqref{thm3} not only for distributions of sums $S=X_{1}+\cdots+X_{n}$ satisfying conditions~\eqref{tau}, but also for distributions from the class $\mathcal{A}_{1}(\tau )$, equivalent to the class of one-dimensional distributions considered by Statulevi\v {c}ius.

Let $\mathcal{A}_{d}(\tau )$, \,$\tau \ge 0$, $d\in \mathbf{N}$,
be the class of $d$-dimensional distributions introduced in an author's paper~\cite{hh29}. The class \,$\mathcal{A}_{d}(\tau )$ \,(with
a fixed \,$\tau \ge 0$) consists of $d$-dimensional
distributions~$F$, \,for which the function
\begin{equation}\label{t1.11}\varphi (z)=\varphi (F,z)=\log \int_{{\mathbf
R}^{d}}e^{\left\langle z,x\right\rangle }F\{dx\}\quad(\varphi
(0)=0)\end{equation}is defined and analytical for \,$ \left\Vert
z\right\Vert \tau <1$, $z\in \mathbf{C}^{d}$, \,and
\begin{equation}\label{b1} \left\vert d_{u}d_{v}^{2}\,\varphi
(z)\right\vert \le \Vert u\Vert \tau \,\bigl\langle\hskip1pt
\mathbb{D}\,v,v\bigr\rangle, \end{equation} for all $u,v\in
{\mathbf R}^{d}$ and $\left\Vert z\right\Vert \tau <1$, where
\,$\mathbb{D}= \text{\rm cov}\,F$ is the covariance operator of distribution~$F$,
 and $ d_{u}\varphi $ is the derivative of  function
$\varphi $ in  direction $u$.

Let's introduce the necessary notation. Below, the symbols $c,c_{1},c_{2},c_{3},\ldots$ will be used for absolute positive constants.
Note that $c$ can be different in different (or even in the same) formulas.
We will write $A\ll B$ if $A\leq c B$. We will also
use the notation $A\asymp B$ if $A\ll B$ and $B\ll A$.
If the corresponding constant depends on, say, $r$, we will write $c(r)$,
$A\ll_r B$, and~$A\asymp_r B$. By $\widehat{F}(t)=\int e^{itx}\,F\{dx\}$, $t\in\mathbf{R}$, we denote the characteristic function of the univariate distribution $F$.

The main result of this paper is contained in the following Theorems \ref{main} and \ref{main1}. They deal with the proximity of univariate distributions.

\begin{theorem}
  \label{main} Let $F=\mathcal L(\xi)\in\mathcal{A}_{1}(\tau )$, $\tau>0$, $\mathbf{E}\,\xi =0$. Then there exists an absolute constant $c_1$ such that
  \begin{equation}\label{thm6}
 W(F,\Phi) =\inf_{\pi}\int\exp\{|x-y|/c_1\tau\}\,d\pi(x,y)\le c_1,
\end{equation}
where $\Phi=\Phi_F$ is the corresponding normal distribution, and the infimum is taken over all two-dimensional probability distributions $\pi$ with marginal distributions $F$ and $\Phi$.
\end{theorem}

\begin{theorem}
  \label{main1} Let $F=\mathcal L(\xi)\in \mathcal{A}_{1}(\tau )$, $\tau>0$, $\mathbf{E}\,\xi =0$. Then there exists an absolute constant $c_2$ such that
\begin{equation}\label{thm33}
 W_\psi(F,\Phi)\le c_2\tau,
\end{equation}
with the Orlicz function $\psi(x)=\exp\{|x|\}-1$,
where $\Phi=\Phi_F$.
\end{theorem}

Theorems \ref{main} and \ref{main1} are equivalent. If $c_1\le2$, then Theorem \ref{main} implies that
\begin{equation}\label{thm68}
\inf_{\pi}\int\exp\{|x-y|/c_1\tau\}\,d\pi(x,y)\le 2.
\end{equation}
If $c_1>2$, then we can choose $c_{3}$
so that $c_1^{c_{3}}=2$, and $c_{3}<1$. Then, by Lyapunov's inequality for moments,
 \begin{equation}\label{thm63}
\inf_{\pi}\int\exp\{c_3\,|x-y|/c_1\tau\}\,d\pi(x,y)\le\inf_{\pi} \Big(\int\exp\{|x-y|/c_1\tau\}\,d\pi(x,y)\Big)^{c_{3}}\le c_1^{c_{3}}=2.
\end{equation}
Now inequalities \eqref{thm68}, \eqref{thm63} imply the statement of Theorem \ref{main1}. It is also obvious that  Theorem~\ref{main1} implies the statement of Theorem \ref{main}.

\section{Properties of classes $\mathcal{A}_{d}(\tau )$}

\label{section 3}

Let's consider the elementary properties of classes $\mathcal{A}_{d}(\tau )$
 (see \cite{hh29, hh31,
hh33,hh34, hh35}). It is easy to see that if
\,$\tau_{1}<\tau_{2}$, \,then \,$\mathcal{A}_{d}\left( {\tau
_{1}}\right) \subset \mathcal{A}_{d}\left( {\tau_{2}}\right) $.
\,Furthermore, the class $\mathcal{A}_{d}\left( {\tau }\right) $
is closed under the convolution operation: if~\,$F_{1},F_{2}\in
\mathcal{A}_{d}\left( {\tau } \right) $, then~$F_{1}F_{2}\, =F_{1}*F_{2}\in \mathcal{A}%
_{d}\left( {\tau }\right) $.  From here on, products and powers of measures are understood
in the sense of convolution.

Let \,$\tau \ge 0$, $F=\mathcal{L}(\xi )\in \mathcal{A}_{d}(\tau )$%
, $y\in \mathbf{R}^{m}$, and $\mathbb{A}:\mathbf{R}^{d}\to \mathbf{R}^{m}$ is a
linear operator. Then
$$
\mathcal{L}(\mathbb{A}\,\xi +y)\in \mathcal{A}_{m}\left( \left\|
\mathbb{A}\right\| \tau \right) ,\quad \mbox{where}\ \ \left\|
\mathbb{A}\right\| =\sup_{x\in \mathbf{R}^{d},\,\left\| x\right\|
\le 1}\left\| \mathbb{A}\,x\right\| .
$$
In particular, for any $a\in \mathbf{R}$
$$
\mathcal{L}(a\,\xi)\in \mathcal{A}_{d}\left( \left|
a
\right| \tau \right) .
$$

The classes \,$\mathcal{A}_{d}(\tau )$ \,are closely related to other
naturally defined classes of multivariate distributions. From
the definition of \,$\mathcal{A}_{d}(\tau )$ it follows that if
\,${\mathcal{L}(\xi )\in \mathcal{A}_{d}(\tau )}$, then
the vector~\,$\xi $ \,has finite exponential moments
\,$\mathbf{E}\,e^{\langle h,\xi \rangle }<\infty $, \,for \,$h\in
\mathbf{R}^{d}$, \,$\left\| h\right\| \tau <1$. \,This leads
to an exponential decay of the distribution tails.

The condition~$\mathcal{L}(\xi )\in \mathcal{A}_{{1}}({\tau )}$
is equivalent to the condition of Statulevi\v cius \cite{hh27}, see also \cite{DJS, 89, 91}, on the growth rate
of the cumulants \,$\gamma_{m}$ \,of a random variable~$\xi$:
\begin{equation}\label{stat}
  |\gamma_{m}|\le \frac{1}{2}\,m!\,\tau ^{m-2}\gamma_{2},\qquad
m=3,4,\dots .
\end{equation}
This equivalence means that if one of these conditions
is satisfied with the parameter~\,$\tau $, then the second is valid with
the parameter~\,$c\tau $, where $c$ is some positive
absolute constant.
Note, however, that the condition \,$\mathcal{L}(\xi )\in \mathcal{A%
}_{{d}}({\tau )}$ differs significantly from other multivariate
analogues of the Statulevi\v cius condition considered by
Rudzkis~\cite{hh23} and Saulis~\cite{hh25}. The review article \cite{91} and the monograph \cite{89}
contain a large number of examples of distributions satisfying  conditions~\eqref{stat} and which are not distributions of sums of a large number of independent terms. Note also that Statulevic\v ius considered more general conditions under which exponential moments are not necessarily finite.

Another class of distributions, denoted by \,$\widetilde{\mathcal
A}_d(\tau)$, \,$\tau \ge 0$, was mentioned in the paper~\cite{hh35}.
It is defined similarly to $\mathcal{A}_{d}(\tau )$ with
replacing~\eqref{b1} with
 \begin{equation}\label{k1.11} \left\vert\,
d_{v}^{2}\,\varphi (z)\right\vert \le 2
\,\bigl\langle\hskip1pt \mathbb{D}\,v,v\bigr\rangle \end{equation}
for all
$v\in {\mathbf R}^{d}$ and $\left\Vert z\right\Vert \tau <1$.
That the classes $\widetilde{\mathcal A}_d(\tau)$ and ${\mathcal A}_d(\tau)$ are also equivalent is easily verified using Cauchy inequalities.
The definition of the classes $\widetilde{\mathcal A}_d(\tau)$ in some sense looks
even more natural than the definition of the classes ${\mathcal A}_d(\tau)$.
The constant 2 in \eqref{k1.11} can be replaced by any other constant $C$, $1<C<\infty$,
bounded away from 1 and  from infinity. The result will also be equivalent classes.

Clearly, the class $\mathcal{A}_{d}\left( {0}\right) $ coincides with the class of all
$d$-dimensional Gaussian distributions. The following inequality~\eqref{pi0} was proved
in the author's paper \cite{hh29} and can be considered as an estimate of the stability of this
characterization:
\begin{equation}
\mbox{if}\ F\in \mathcal{A}_{d}(\tau ),\ \mbox{then}\ \pi \left( F,\,\Phi
_F\right) \le c\,d^{2}\tau \,\log ^{*}(\tau ^{-1}),  \label{pi0}
\end{equation}
where \,$\pi (\,\cdot\,,\,\cdot\,)$ is the Prokhorov distance,
defined in \cite{hh19}, and $ \Phi _F $ denotes
the Gaussian distribution whose mean and covariance operator
are the same as those of~\,$ F $.
Here $\log ^\ast b=\max \left\{
1,\log b\right\} $ for $b>0$, and $\log$ is used to
denote the natural logarithm.
Note that Theorems~\ref{main} and \ref{main1} of this paper
can also be viewed as stability estimates for the above-mentioned
characterization of Gaussian distributions in transport metrics, and so far in the one-dimen\-sional case.

The Prokhorov distance between
distributions $ F, G $ can be determined by the formula
\[
\pi (F,G)=\inf \left\{ \lambda :\pi (F,G,\lambda )\leq \lambda \right\} ,
\]
where
\[
\pi (F,G,\lambda )=\sup_{Y}\max \left\{ F\{Y\}-G\{Y^{\lambda
}\},G\{Y\}-F\{Y^{\lambda }\}\right\} ,\quad \lambda >0,
\]
and $Y^{\lambda }=\{y\in \mathbf{R}^{d}:\inf\limits_{x\in Y}\left\|
x-y\right\| <\lambda \}$ is the $ \lambda $-neighborhood of a Borel
set~$ Y $\break (see~{\cite{BMS, 5.}}).

In the author's paper \cite{hh29} it was also established that
\begin{equation}
\mbox{if}\ F\in \mathcal{A}_{d}(\tau ),\ \mbox{then}\ \pi (F,\Phi
_F,\lambda )\le c\,d^{2}\exp \Big\{-\frac{\lambda }{c\,d^{2}\tau
}\Big\},\quad \lambda >0.  \label{pi}
\end{equation}
What is important here is that  inequality~\eqref{pi} is proved for
all~\,${\tau >0}$ \,and for an arbitrary covariance operator~cov$\,F$.

By the Strassen--Dudley theorem (see Dudley \cite{hh9}) and
according to  inequality~\eqref{pi}, for any distribution $F\in
\mathcal{A}_{d}(\tau )$ and any $\lambda>0$, one can construct random vectors \,$\xi $ \,and \,$\eta $ \,on the same probability
space with
\,$\mathcal{L}(\xi )=F$ \,and \,$\mathcal{L}(\eta )=\Phi_F$, \,so
that
\begin{equation}
\mathbf{P}\left\{ \Vert \xi -\eta \Vert >\lambda \right\}=\pi
(F,\Phi_F,\lambda ) \le c\,d^{2}\exp %
\Big\{-\frac{\lambda }{c\,d^{2}\tau }\Big\}.  \label{pi1}
\end{equation}
We emphasize that the Strassen--Dudley theorem guarantees the existence
of a construction with equality in \eqref{pi1} only for a fixed $\lambda$. An example showing  impossibility
of a construction with equality in \eqref{pi1} for all $\lambda$
simultaneously can be found in the survey~\cite{BMS}. The Strassen--Dudley theorem enables us to automatically derive statements
of the type~\eqref{pi1} from estimates for $\pi (F,G,\lambda)$. Strassen's original proof \cite{Str}
was non-constructive. Dudley~\cite{hh9} gave a compli\-ca\-ted constructive
proof based on combinatorial ideas. Finally, Schay
\cite{Schay} found a short proof relying on the duality theorem.

If  equality \eqref{pi1} were proven for all $\lambda>0$ simultaneously on the same probability space, then the assertion of Theorem~\ref{main} would automatically follow from it,
for any dimension~$d$, $1\le d<\infty$. Therefore,  inequality \eqref{pi} gives grounds to expect a possibility to generalize Theorems~\ref{main}
and~\ref{main1} to the multidimensional case.

If \,$F$ is an infinitely divisible distribution with
spectral measure concentrated on the ball $\left\{ x\in \mathbf{R}%
^{d}:\left\| x\right\| \le \tau \right\} $, then ${F\in
\mathcal{A}_{d}(c\tau )}$, where $c$ is some positive
absolute constant. In the paper \cite{hh29}, one can find less
restrictive conditions on  spectral measure, ensuring
that an infinitely divisible distribution
belongs to  ~${\mathcal{A}_{d}(c\tau )}$.

In particular, for the Poisson distribution $\Pi_{\lambda}$ with parameter $\lambda>0$, the inclusion $\Pi_{\lambda}\in \mathcal{A}_{1}(c )$ holds.
It follows from Theorem~\ref{main1} that
\begin{equation}\label{thm44}
\sup_{\lambda} W_\psi(\Pi_{\lambda},\Phi_{\Pi_{\lambda}})\le c,\quad\hbox{for }\psi(x)=\exp\{|x|\}-1.
\end{equation}
This is the statement of Corollary 2.2 of Rio~\cite{Rio}. But Theorem~\ref{main1} contains a more general assertion. In \eqref{thm44} we can replace the set of Poisson distributions by the set of all infinitely divisible distributions with the L\'evy--Khinchin spectral measures concentrated on the interval~$[-1,1]$.

Distributions from the classes $\mathcal{A}_{d}(\tau )$ are directly used in the formulations of the author's results
\cite{hh33}--\cite{hh35} on estimating the accuracy of strong Gaussian approximations for sums of independent random vectors
in the most important case, when
the summands have finite exponential moments (see also the review article \cite{zai13}). Multivariate analogues of the one-dimensional results of
Sakhanenko~\cite{hh24}, who generalized and significantly refined
the results of Koml\'os, Major, and Tusn\'ady~\cite{hh17}, for the case of non-identically distributed random variables, were obtained. Sakhanenko considered the following classes of univariate distributions:
\[
\mathcal{S}_{1}(\tau )=\left\{ \mathcal{L}(\xi ):\ \mathbf{E}\,\xi =0,\
\mathbf{E}\,|\xi |^{3}\exp \left\{|\xi |/\tau \right\} \le \tau \,\mathbf{%
E}\,|\xi |^{2}\right\}, \quad\tau>0.
\]
In the author's preprint~\cite{Zai84} it was noted that the classes
$\mathcal{S}_{1}(\tau )$ are equivalent to the classes of distributions
$\mathcal{B}_{1}(\tau )$ satisfying the conditions of the
Bernstein inequality~(see definition~\eqref{bern}), in the sense
that if one of the conditions for membership in a class is satisfied with
a parameter~\,$\tau $, then the second is true with
 parameter~\,$c\tau $, where $c$ is some positive
absolute constant. Sakhanenko's results~\cite{hh24} were formulated as estimates of  exponential moments of maximal deviation
of sums of independent random variables with distributions in $\mathcal{S}_{1}(\tau )$, constructed on a  probability space, from  corresponding sums of independent normally distributed terms. The form of Sakhanenko's estimates is almost the same as in~\eqref{thm6}, so that they imply an analogue of inequality~\eqref{thm6} with  right-hand side replaced by $c(1+\sigma/\tau)$ for the distributions of sums of independent random variables with distributions from $\mathcal{S}_{1}(\tau )$. Here, $\sigma^{2}$ denotes the variance of the sum under consideration. For convolutions of distributions from $\mathcal{S}_{1}(\tau )$, some estimates of  moments of exponential type in the central limit theorem are contained in Sakhanenko~\cite{sah96}.

In  the author's papers~\cite{Zai84} and \cite{hh30},
 inequalities \eqref{pi0} and \eqref{pi} (with ${d}^{2}$ replaced by
${d}^{5/2}$) were proved for convolutions of distributions from the class
$\mathcal{B}_{d}(\tau )$, where ${\tau >0}$ and
\begin{eqnarray}
\mathcal{B}_{d}(\tau ) &=&\Big\{F=\mathcal{L}(\xi
):\mathbf{E}\,\xi =0,\ \left| \mathbf{E}\,\left\langle \xi
,v\right\rangle ^{2}\left\langle \xi
,u\right\rangle ^{m-2}\right| \label{bern}\\
&\le &\frac{1}{2}m!\,\tau ^{m-2}\left\| u\right\| ^{m-2}\,\mathbf{E}%
\,\left\langle \xi ,v\right\rangle ^{2}\ \ \mbox{for all}\ u,v\in \mathbf{R}%
^{d},\ m=3,4,\dots \Big\},\nonumber
\end{eqnarray}
satisfying multidimensional analogues of the Bernstein inequality conditions. Sa\-kha\-nen\-ko's condition $\mathcal{L}(\xi )\in
\mathcal{S}_{{1}}({\tau )}$ is equivalent to the condition $\mathcal{L}(\xi
)\in \mathcal{B}_{{1}}({\tau )}$. Note that\break if $F\left\{
\left\{ x\in \mathbf{R}^{d}:\left\| x\right\| \le \tau \right\}
\right\} =1$, $\mathbf{E}\,\xi =0$, then $F\in \mathcal{B}_{d}(\tau )$.

Let us formulate  relations between  classes \,$\mathcal{A}_{d}(\tau )$
\,and $\mathcal{B}_{d}(\tau )$. Let $\sigma ^{2}_F$
denote the maximal eigenvalue of the covariance operator
of  distribution  \,$F $. Then \medbreak

\noindent a) If $F=\mathcal{L}(\xi )\in \mathcal{B}_{d}(\tau )$,
then $\sigma ^{2}_F\le 12\,\tau ^{2}$, $\mathbf{E}\,\xi =0$ and ${F\in \mathcal{A}%
_{d}(c\tau )} $.

\noindent b) If $F=\mathcal{L}(\xi )\in \mathcal{A}_{d}(\tau )$,
$\sigma
^{2}_F\le \tau ^{2}$ and $\mathbf{E}\,\xi =0$, then $F\in \mathcal{B}%
_{d}(c\tau )$. {\medbreak}

In particular, the distribution of the sum $S=X_{1}+\cdots+X_{n}$ under  conditions \eqref{tau} belongs to the class ${ \mathcal{A}%
_{d}(c\tau )} $.

Thus, roughly speaking, $ \mathcal{B}_{d}(\tau )$
forms a subclass of distributions $F=\mathcal{L}(\xi-
\mathbf{E}\,\xi)$ such that ${\mathcal{L}(\xi )\in \mathcal{A}_{d}(c\tau
)}$, and $\sigma ^{2}_F\le 12\,\tau ^{2}$. Inequalities \eqref{pi0} and \eqref{pi}
in this case indicate only that both
distributions being compared are close to  degenerate law~$E$ concentrated at the origin. If
${F=\mathcal{L}(\xi )\in \mathcal{A}_{d}(\tau )}$ and $\sigma
^{2}_F$ is significantly greater than $\tau ^{2}$, then
${\mathcal{L}(\xi/\sigma_F )\in \mathcal{A}_{d}(\tau/\sigma_F
)}$ and  inequalities \eqref{pi0} and \eqref{pi} reflect the proximity
of  distribution  $F$ to the corresponding Gaussian law.

Let \;$\tau\ge0$, \,$F=\mathcal L(\xi)\in\mathcal A_d(\tau)$, \,
$\| h\|\tau<1$, \,$h\in {\mathbf R}^d$. \;
We define the distribution \,
$\overline F=\overline F(h)$ \,by  relation
$$
\overline F\{dx\}=
\Bigl(\!\mathbf E\, e^{\langle h,\xi\rangle}\!\Bigr)^{-1}e^{\langle h,x\rangle }F\{dx\}.
$$
We denote by $\overline\xi=\overline\xi(h)$ a random vector with  distribution $\mathcal L(\overline\xi(h))=\overline F(h)$.
Distributions $\overline F(h)$ are sometimes called Cram\'er transforms (or Esscher transforms, see \cite{BG}). In
the proofs in  \cite{hh29}, \cite{hh30}, {\cite{hh33}--\cite{hh35}},  distributions
$\overline F(h)$ are used to estimate  probabilities of large deviations,
corresponding to the conditional densities. Another important property
of  classes $\mathcal{A}_{d}(\tau )$ is that ${\overline F(h)\in
\mathcal{A}_{d}(2\tau )}$ for $\left\| h\right\|\tau \le 1/2$, see item b) of Lemma \ref{l1}. This
makes it possible to systematically apply the results obtained
for the original distributions to their Cram\'er transforms and thereby refine the estimates.

 Kolmogorov \cite{K} posed the problem of estimating the accuracy of  infinitely divisible approxi\-ma\-tion of  distributions of sums of independent random variables whose distributions are con\-cen\-trated on short intervals
of small length $\tau\le1/2 $  up to a small probability~$p$. In the particular case when $p=0$, we are talking about approximating the distributions of sums $S=X_{1}+\cdots+X_{n}$ for $d=1$ and under  conditions \eqref{tau}.
In this case, Kolmogorov \cite{K,k63} obtained the estimate
\begin{equation}
L(F,\Phi_F)\ll
\tau^{1/2}\log^{1/4}(1/\tau),
\label{infdiv2}
\end{equation}
where $L(\, \cdot \,,\,\cdot \,)$ is the L\'evy distance. From the above it follows that $F\in\mathcal{A}
_{1}(c\tau ) $ and  inequalities \eqref{pi0} and \eqref{pi} can be viewed as  generalizations, and refinements of  inequality~\eqref{infdiv2}.
Note that the formulations  in  \cite{K,k63} differ from~\eqref{infdiv2}. To deduce this inequality from them, an elementary additional analysis is required.

Let $X,X_1,\ldots,X_n$ be independent identically distributed random variables such that
\begin{equation}\label{tau4}
\mathbf{E}\,X=0\quad\hbox{and}\quad\mathbf{E}\,\exp\{t|X|\}<\infty\quad\hbox{for some }\quad t>0.
\end{equation}
Then it is easy to verify that there exists $c(F)$ such that $F=\mathcal{L}(X)\in\mathcal{A}
_{1}(c(F))$,
and  distribution of the normalized sum $F_{n}=\mathcal{L}\big((X_{1}+\cdots+X_{n})/\sqrt{n}\big)$ belongs to $\mathcal{A}
_{1}(c(F)/\sqrt{n})$. Applying Theorem~\ref{main1}, we obtain that
\begin{equation}\label{thm4}
 W_\psi(F_n,\Phi_{F_n})\ll_F1/\sqrt{n},
\end{equation}
with  Orlicz function $\psi(x)=\exp\{|x|\}-1$. This is the statement of Theorem 2.1 of Rio~\cite{Rio}. We emphasize once again that in the main results of this article we consider univariate distributions that satisfy condition \eqref{b1} on the Laplace transforms and, generally speaking, cannot be represented as convolutions of a large number of identical distributions. In this case,  condition \eqref{b1} turns into
\begin{equation}\label{b15} \left\vert \,\varphi'''
(z)\right\vert \le  \tau \,\sigma^{2}\quad\text{for }|z|\tau\le1. \end{equation}

From the above it follows that the assertions of Theorems \ref{main} and \ref{main1} are also  valid for convolutions of univariate distributions concentrated on the interval $[-\tau,\tau]$ or satisfying the conditions of Bernstein's inequality, as well as for infinitely divisible distributions with their L\'evy--Khinchin spectral measures concentrated on the same interval, and for distributions satisfying the Statulevi\v cius conditions~\eqref{stat}.
In terms of content and methods of proof, they can be viewed as simply and clearly formulated statements from large deviation theory.

\section{Proof of Theorem \ref{main}}

We will need the following Lemma~\ref{l1} about the properties of the Cram\'er transform (it is contained in~\cite [Lemmas 2.1, 3.1]{hh29}).

 \begin{lemma}\label{l1}
Let \;$\tau\ge0$, \;
$F=\mathcal L(\xi)\in\mathcal A_d(\tau)$, \,$\mathbf E\, \xi=0$, \;$h\in {\mathbf R}^d$, \,${\|h\|\tau<1}$,
 \, $\mathbb D=\hbox{\rm cov} F$,
\,$ \mathbb D(h)=\hbox{\rm cov} {\overline F}(h)$. \,Denote by $\sigma^{2} $ the minimal eigenvalue of $\mathbb D$. Then

a$)$ for any \,$u\in {\mathbf R}^d$ \,the following relations hold$:$

\begin{equation}\label{1}
\langle\mathbb D(h)\,{u},\,{u}\rangle=\langle\mathbb D\,{u},\,{u}\rangle\bigl(1+\theta\|h\|\tau\bigr),
\end{equation}
\begin{equation}
\label{5}\log \mathbf E\,e^{i\langle h,\xi\rangle}=-\frac12\,\langle\mathbb D\,{h},\,{h}\rangle\Bigl(1+\frac13\,\theta\|h\|\tau\Bigr),
\end{equation}\begin{equation}
\label{6}\log \mathbf E\,e^{\langle h,\xi\rangle}=\frac12\,\langle\mathbb D\,{h},\,{h}\rangle\Bigl(1+\frac13\,\theta\|h\|\tau\Bigr)
\end{equation}
$($here and below \;$\theta$ \,
symbolizes various quantities not exceeding one in absolute value$:|\theta|\le1);$

b$)$ If \;$\|h\|\tau\le1/2$, \;then \;
$\overline F(h)\in\mathcal A_d(2\tau)$;

c$)$ For $x\in\Pi=\big\{x\in  {\mathbf R}^d:4.8\,\tau\sigma^{-1}\,\bigl\|\mathbb{D}^{-1/2}x\bigr\|\le1\big\}$
there exists a parameter $h=h(x)\in  {\mathbf R}^d$ such that

 \begin{equation}
\mathbf E\, \overline\xi(h)=x,
\end{equation}
\begin{equation}
\|h\|\tau\le1/2,
\end{equation}
\begin{equation}
\sigma\,\|h\|\le\bigl\|\mathbb{D}^{1/2}h\bigr\|\le 2.4\,\bigl\|\mathbb{D}^{-1/2}x\bigr\|,
\end{equation}
\begin{equation}
\bigl\|\mathbb{D}^{1/2}h-\mathbb{D}^{-1/2}x\bigr\|\le
2.88\,\theta\tau\sigma^{-1}\,\bigl\|\mathbb{D}^{-1/2}x\bigr\|^2,
\end{equation}
\begin{equation}
\mathbf E\,\exp\big\{\langle h,\xi\rangle-\langle h,x\rangle\big\}=\exp\Big\{-\frac{\,1\,}{2}\,\bigl\|\mathbb{D}^{-1/2}x\bigr\|^2+
10.08\,\theta\tau\sigma^{-1}\,\bigl\|\mathbb{D}^{-1/2}x\bigr\|^3\Bigr\}.
\end{equation}
\end{lemma}

Below
\begin{equation}
\Xi (x)=e^{x^2/2}\int_{x}^{\infty}e^{-y^2/2}\,dy, \quad x>0,
\end{equation} is the Mills ratio.
We will need the following lemma (see \cite[Lemma 1.2 of Chapter VI]{Arak and Zaitsev}).

\begin{lemma}\label{l3}
  Let $x,\varepsilon>0$. Then
\begin{equation}\label{554}
  0\le\Xi (x)-\Xi (x+\varepsilon)\le\frac{\varepsilon}{x^{2}},
\end{equation}
\begin{equation}\label{556}
\Xi (x)=\frac{1}{x}\Big(1-\frac{|\theta|}{x^{2}}\Big).
\end{equation}
Let $\Phi_{\sigma}$ be the univariate Gaussian distribution with zero mean and variance~$\sigma^{2}$. Then
  \begin{equation}\label{527}
1- \Phi_\sigma(x+\varepsilon)\le (1- \Phi_\sigma(x))\,\exp\Big\{-\frac{2x\varepsilon+\varepsilon^2}{2\sigma^{2}}\Big\}.
\end{equation}Hence, for $x>\varepsilon$
\begin{equation}\label{5279}
1- \Phi_\sigma(x)\le (1- \Phi_\sigma(x-\varepsilon ))\,\exp\Big\{-\frac{2x\varepsilon-\varepsilon^2}{2\sigma^{2}}\Big\}.
\end{equation}
\end{lemma}

Let $\rho(F, \Phi)=\sup_{x}\bigl|F(x)-\Phi(x)\bigr|$ be the Kolmogorov distance, uniform distance between distribution functions.
\begin{lemma}\label{l2}
Under the conditions of Theorem $\ref{main}$,
\begin{equation}\label{44}
  \rho(F, \Phi)\ll\tau/\sigma,
\end{equation}
where $\sigma^{2}=\text{\rm Var}\, \xi$ is the common variance of distributions $F$ and $\Phi$.
\end{lemma}

\noindent{\bf Proof.}
Using inequality
\begin{equation}\label{2}
  \bigl|e^{z_1}-e^{z_2}\bigr|\le\bigl|{z_1}-{z_2}\bigr|\,\max\big\{\bigl|e^{z_1}\bigr|,\bigl|e^{z_2}\bigr|\big\},
\end{equation}
valid for $z_1, z_2\in\mathbb{C}$, and applying inequality \eqref{5}, we obtain that for $|t|\tau\le1$
\begin{equation}\label{3}
  \bigl|\widehat{F}(t)-\widehat{\Phi}(t)\bigr|\le\frac\tau6\sigma^2|t|^3\exp\Big\{-\frac13\sigma^2t^2\Big\}.
\end{equation}

Therefore, using the standard smoothing inequality (see \cite[Theorem 1.2 of Chapter~III]{Arak and Zaitsev}), we find that for $T=1/\tau$
\begin{equation}
\rho(F, \Phi)\ll\int_{0}^{T}\biggl|\frac{\widehat{F}(t)-\widehat{\Phi}(t)}{t}\biggr|\,dt+\frac{1}{\sigma T}\ll\frac\tau\sigma.
\end{equation}

\def\t{\tau}The following lemma contains an analogue of Bernstein's inequality for distributions from the class $\mathcal{A}_{1}(\tau )$.

\begin{lemma}\label{2.6}
Let, under the conditions of Theorem $\ref{main}$, \,$\text{\rm Var}\, \xi=\sigma^2$. \, Then
\begin{equation}\label{be}
{\mathbf P}\bigl\{\xi
\ge x\bigr\}\le\max\Bigl\{\exp\Bigl\{-\frac{x^2\!}{4
\,\sigma^2}\Bigr\},\;
\exp\Bigl\{-\frac{x}{4\,\t}\Bigr\}\Bigr\},\qquad x\ge0.
\end{equation}
\end{lemma}
The proof of this lemma almost literally repeats the proof of Bernstein's inequality. Let $0\le h\tau\leqslant\frac{1}{2}$.
By~\eqref{6},
\begin{equation}\label{00}
{\mathbf E}\,e^{h\xi}\leqslant  \exp \{h^2\sigma^2\}
\nonumber\end{equation}
 and
\begin{equation}
{\mathbf P}\{\xi\geqslant x\}\leqslant e^{-hx}\,{\mathbf E}\,e^{h\xi}\leqslant  \exp \{h^2\sigma^2-hx\}.
\nonumber\end{equation}
Let's choose the parameter $h$ depending on $x$. If $0\leqslant x\leqslant\frac{\sigma^2}{\tau}$, we take $h=\frac{x}{2\sigma^2}$ and obtain the bound
\begin{equation}\label{be1}
{\mathbf P}(S\geqslant x)\leqslant \exp \left\{-\frac{x^2}{4\,\sigma^2}\right\}.
\end{equation}
And if $x>\frac{\sigma^2}{\tau}$, take $h=\frac{1}{2\tau}$ and get
\begin{equation}\label{be2}
{\mathbf P}(S\geqslant x)\leqslant \exp \left\{\frac{\sigma^2}{4\tau^2}-\frac{x}{2\tau}\right\}\leqslant \exp \left\{-\frac{x}{4\,\tau}\right\}
.\end{equation}
 Now the inequality \eqref{be} follows from \eqref{be1} and \eqref{be2}.\medskip

\noindent{\bf Proof of Theorem \ref{main}.} Without loss of generality, we assume that the distribution function $F$ is infinitely differentiable and strictly increasing.
To justify this, it suffices to consider, instead of the distribution $F$, the convolution of this distribution with a Gaussian distribution having zero mean and
positive variance tending to zero, and also to use the standard tool for proving theorems on strong
approximation, Lemma A from Berkes and Philipp \cite{hh3}; see, for example, the proof of Theorem 3.1 in Rio \cite{Rio}. Thus, under these assumptions, the strictly increasing inverse function $F^{-1}(\,\cdot\,)$ is well defined.

Let a random variable $\eta$ have the distribution \,$\mathcal L(\eta)=\Phi$. Write ${\xi=F^{-1}\big(\Phi(\eta)\big)}$. It is clear that \,$\mathcal L(\xi)=F$
and ${\eta=\Phi^{-1}\big(\Phi(\eta)\big)}$. This means that the random variables \;$\xi$ \;and \;$\eta$ \;are defined as the Smirnov transforms of a random variable $\Phi(\eta)$ uniformly distributed on the interval~\emph{}$[0,1]$. This is exactly how random variables with given distributions are constructed in the proof of  equality \eqref{val} in \cite{V}. Then if the random variable $\xi$ takes  some specific value~$x\in\mathbf{R}$, the random variable $\eta$ will take the value $\Phi^{-1}\big(F(x)\big)$.

Further reasoning is carried out under the assumption that $\xi=x$ and $\tau\le c_4\sigma$, where $\sigma^2=\text{\rm Var}\, \xi$, and the choice of $c_4$ will be refined during the proof.

First we consider the case when $|x|\le2\sigma$. Let
\begin{equation}\label{50gg}
\phi(u)=\frac{1}{\sqrt{2\pi}\,\sigma}\,\exp\Bigl(-\frac{u^2}{2
\,\sigma^2}\Bigr)
\end{equation}
be the density of distribution $\Phi$. Recall that according to Lemma \ref{l2}
\begin{equation}\label{49}
\rho(F, \Phi)\le c_5\frac\tau\sigma.
\end{equation}
Let $|u|\le2\sigma$, $|y|\le3\sigma$. It is obvious that then
\begin{equation}\label{50ppp}
\bigl| \Phi(u)- \Phi(y)\bigr|\ge|u-y|\,\phi(3\sigma)=\frac{e^{-9/2}\,|u-y|}{\sqrt{2\pi}\,\sigma}=\frac{c_6\,|u-y|}{\sigma}, \quad c_6=\frac{e^{-9/2}}{\sqrt{2\pi}},
\end{equation}
\begin{equation}\label{51}
\bigl| F(u)-F(y)\bigr|\ge\frac{c_6\,|u-y|}{\sigma}-2c_5\frac\tau\sigma.
\end{equation}

\begin{lemma}\label{quant}
There exist absolute positive constants $c_4,c_7$ such that for $\tau\le c_4\sigma$, $|x|\le2\sigma$
the following inequalities hold:
 \begin{equation}\label{52}
 \Phi(x+c_7\tau)\ge F(x),\quad  F(x+c_7\tau)\ge \Phi(x).
\end{equation}
\end{lemma}

Indeed, according to \eqref{49}--\eqref{51}, for $|x|\le2\sigma$
 \begin{equation}
 \Phi(x+c_7\tau)- F(x)\ge F(x+c_7\tau)- F(x)-c_5\frac\tau\sigma\ge\frac{c_6\,c_7\tau}{\sigma}-3c_5\frac\tau\sigma\ge0,
\end{equation}
if we choose $c_7=3c_5/c_6$ and if $x+c_7\tau\le3\sigma$. The last inequality becomes obvious if $c_7\tau\le\sigma$. For this, it is sufficient to choose $c_4\leq c_7^{-1}$. The second inequality in \eqref{52} is verified similarly.

Thus, according to Lemma \ref{quant},
 \begin{equation}\label{53}
|\,\xi-\eta\,|< c_7\tau\emph{},\qquad \text{if}
\quad  |\xi|\le 2\sigma.
 \end{equation}

Let $2\sigma\le x\le\sigma^2/5\tau$, and the parameter $h=h(x)\in {\mathbf R}$ be chosen in accordance with item~c)
of the one-dimensional version
of Lemma~\ref{l1} (whose condition
$x\in\Pi$ is satisfied) and it is such that
$\mathbf{E}\,\overline \xi(h)=x$,  $\|h\|\tau\le1/2$, $\mathcal L(\overline\xi(h))=\overline F=\overline F(h)$,
\begin{equation}\label{50}
\left|\sigma h-x/\sigma\right|\le2.88\,\tau\sigma^{-1}x^2\sigma^{-2}.
\end{equation}
According to  relation \eqref{1} of Lemma~\ref{l1},
\begin{equation}\label{99}\sigma^2(h)=
\text{\rm Var}\,\overline \xi(h)=\sigma^2\bigl(1+\theta\|h\|\tau\bigr).
\end{equation}

Introduce the distribution ${H}=\Phi_{\overline{F}}$. Then, in accordance with item~b) of Lemma~\ref{l1}, $\overline F(h)\in\mathcal A_d(2\tau)$ and,
by Lemma~\ref{l2} taking into account \eqref{99},
\begin{equation}\label{98}
\rho({H},\overline{F})\ll\frac\tau{\sigma(h)}\ll\frac\tau{\sigma}.
\end{equation}
Integrating by parts, we obtain
\begin{equation}\label{7}
1-F(x)=\mathbf{E}\,e^{h\xi}\int_{x}^{\infty}e^{-hy}\overline F\{dy\}=
\mathbf{E}\,e^{h\xi}\bigg(\int_{x}^{\infty}h\,e^{-hy}\overline F(y)\,dy-e^{-hx}\,\overline F(x)\bigg),
\end{equation}
\begin{equation}
\int_{x}^{\infty}e^{-hy}H\{dy\}=
\int_{x}^{\infty}h\,e^{-hy}H(y)\,dy-e^{-hx}\,H(x).
\end{equation}
On the other hand, it is easy to verify that
\begin{equation}\label{8}
\int_{x}^{\infty}e^{-hy}H\{dy\}=
\frac{1}{\sqrt{2\pi}}e^{-hx}\,\Xi(h\sigma(h)).
\end{equation}

From \eqref{98}--\eqref{8} it follows that
\begin{multline}\label{9}\biggl|\int_{x}^{\infty}e^{-hy}\overline F\{dy\}-\frac{1}{\sqrt{2\pi}}e^{-hx}\,\Xi(h\sigma(h))\biggr|
   \\ \le\biggl|\int_{x}^{\infty}h\,e^{-hy}(\overline F(y)-  H(y))\,dy-e^{-hx}(\overline F(x)-  H(x))\biggr|\\
   \le2\,e^{-hx}\,\rho({H},\overline{F})\ll e^{-hx}\,\frac\tau\sigma.
   \end{multline}
Applying  inequality \eqref{554} of Lemma \ref{l3}, we obtain
   \begin{equation}
\biggl|\Xi(h\sigma(h))-\Xi(h\sigma)\biggr|\ll\frac{h|\sigma(h)-\sigma|}{h^2\sigma^2}\ll\frac\tau\sigma.
\end{equation}

If $x\ge2\sigma$, then $h\sigma\ll x/\sigma$ and
$\Xi(h\sigma)\gg\Xi(x/\sigma)\gg\sigma/x$.  Applying inequality \eqref{554} again, as well as \eqref{50}, we obtain
\begin{equation}
\bigl|\Xi(h\sigma)-\Xi(x/\sigma)\bigr|\ll\frac{\bigl|h\sigma-x/\sigma\bigr|\,\sigma^2}{x^2}\ll\frac\tau\sigma.
\end{equation}
Hence,
 \begin{equation}
\biggl|\int_{x}^{\infty}e^{-hy}\overline F\{dy\}-\frac{1}{\sqrt{2\pi}}e^{-hx}\,\Xi(x/\sigma)\biggr|\ll\frac\tau\sigma\,e^{-hx}.
\end{equation}

Applying the  above inequalities, we obtain that
\begin{multline}\label{10}
  1-F(x)=\mathbf{E}\,e^{h\xi}\bigg(\frac{1}{\sqrt{2\pi}}e^{-hx}\,\Xi(x/\sigma)+\int_{x}^{\infty}e^{-hy}\overline F\{dy\}-\frac{1}{\sqrt{2\pi}}e^{-hx}\,\Xi(x/\sigma)\bigg) \\ =\frac{1}{\sqrt{2\pi}}\mathbf{E}\,e^{h\xi-hx}\,\Xi(x/\sigma)\Big(1+\theta\,c\,\frac\tau\sigma\,\frac x\sigma\Big) \\ =(1-\Phi(x))\exp\Big\{\theta\,c_{8}\,\frac\tau\sigma\,\frac {x^3}{\sigma^3}\Big\} \\ = \frac{1}{\sqrt{2\pi}}\,e^{-x^2/2\sigma^2}\,\Xi(x/\sigma)\exp\Big\{\theta\,c_{8}\,\frac\tau\sigma\,\frac {x^3}{\sigma^3}\Big\}
  .
\end{multline}

\begin{lemma}\label{l6}
There exist absolute positive constants $c_9,\ldots,c_{11}$ such that
\begin{equation}\label{55}
  1-\Phi(x+\beta(x))\le1-F(x)\le
 1-\Phi(x-\beta(x))
\end{equation}
 for $\tau/\sigma\le c_9$, $2\sigma\le x\le z= c_{10}\sigma^2/\tau$, where
 \begin{equation}\label{56}
 \beta(x)=c_{11}\tau\,x^2\sigma^{-2}.
\end{equation}
\end{lemma}

We set $c_{11}=4c_{8}$. Then, by choosing a sufficiently small $c_{10}$, we ensure that the inequality $\beta(x)\le x$ is satisfied. Now applying  inequalities
\eqref{527} and \eqref{5279}
with $\varepsilon=\beta(x)$, we obtain
\begin{multline}\label{a90}
  1-F(x)=(1-\Phi(x))\exp\Big\{\theta\,c_{8}\,\frac\tau\sigma\,\frac {x^3}{\sigma^3}\Big\}\\
  \le (1- \Phi(x-\beta(x) ))\exp\Big\{-\frac{(2x-\beta(x))\beta(x)}{2\sigma^{2}}+c_{8}\,\frac\tau\sigma\,\frac {x^3}{\sigma^3}\Big\}
  \le 1- \Phi(x-\beta(x) ),
\end{multline}
\begin{multline}\label{a904}
  1-F(x)=(1-\Phi(x))\exp\Big\{\theta\,c_{8}\,\frac\tau\sigma\,\frac {x^3}{\sigma^3}\Big\}\\
  \ge (1- \Phi(x+\beta(x) ))\exp\Big\{\frac{(2x+\beta(x))\beta(x)}{2\sigma^{2}}-c_{8}\,\frac\tau\sigma\,\frac {x^3}{\sigma^3}\Big\}\ge 1- \Phi(x+\beta(x) ),
\end{multline}
completing the proof of the lemma.

Applying Lemma~\ref{l6}, we obtain that
\begin{equation}\label{555}
|\,\xi-\eta\,|< c_{11}\tau\xi^2/\sigma^2,\qquad \text{if}
\quad 2\sigma\le \xi\le z= c_{10}\sigma^2/\tau, \quad \tau/\sigma\le c_9.
\end{equation}

To prove Theorem~\ref{main} it suffices to prove
that the absolute constant $c_{14}$ can be chosen so large that ${\mathbf E}\,\exp\big\{|\,\xi-\eta\,|/c_{14}\tau\big\}\ll1$.

First, we will assume that $\tau/\sigma\le c_9$.
It is clear that
\begin{multline}\exp\big\{|\,\xi-\eta\,|/c_{14}\tau\big\}
\le \exp\big\{|\,\xi-\eta\,|/c_{14}\tau\big\}\,
{\bf1} \big\{ |\xi|\le 2\sigma\big\}\\
+\exp\big\{|\,\xi-\eta\,|/c_{14}\tau\big\}\,
{\bf1} \big\{2\sigma\le |\,\xi\,|\le z\big\}\\
+\exp\big\{|\,\eta\,|/c_{14}\tau+|\,\xi\,|/c_{14}\tau\big\}\,{\bf1} \big\{ |\xi|\ge z\big\}.
\end{multline}
According to \eqref{53}, for $c_{14}>c_7$
$$
{\mathbf E}\,\exp\big\{|\,\xi-\eta\,|/c_{14}\tau\big\}\,
{\bf1} \big\{ |\xi|\le 2\sigma\big\}\le e.
$$

By Lemma~\ref{2.6},
\begin{equation}\label{be44}
{\mathbf P}\bigl\{\xi^{2}
\ge x^{2}\bigr\}={\mathbf P}\bigl\{|\,\xi\,|
\ge x\bigr\}\le2\,\max\Bigl\{\exp\Bigl\{-\frac{x^2\!}{4
\,\sigma^2}\Bigr\},\;
\exp\Bigl\{-\frac{x^{2}}{4\,c_{10}\,\sigma^{2}}\Bigr\}\Bigr\}
\end{equation}for \;$ 0\le x\le z= c_{10}\,\sigma^{2}/\tau$.
Set
$$
W=
c_{11}\xi^2/\sigma^2.
$$

Using \eqref{be44}, we obtain that there exists $c_{13}$ such that
\begin{eqnarray*}
\mathbf P\bigl\{W\ge u\bigr\}
&\le&
\mathbf P\Bigl\{\bigl|\xi\bigr|
\ge \sigma\,\sqrt {u/c_{11}}\Bigr\}\\
&\le&2\,\exp\Bigl\{- u/c_{13}\Bigr\}
\end{eqnarray*}
for \;$0\le u\le\gamma=c_{11} z^2/\sigma^{2}$. Set $c_{12}=2c_{13}$, $v=1/c_{12}$.
Integrating by parts, we obtain
\begin{eqnarray*}
{\mathbf E\exp\big\{v\,W}\big\}\,{\bf1} \big\{ |\,\xi\,|\le z\big\}&\le&1+\int_{0}^{\gamma}
v\,e^{v\,u}\,\mathbf P\bigl\{W\ge u\bigr\}\,du\\
&\le&1+\frac2{c_{12}}\int_0^{ \gamma}
e^{u/c_{12}}\,\exp\Bigl\{- u/c_{13}\Bigr\}\,{du}\\
&\le&1+\frac2{c_{12}}\int_0^{\infty}
e^{-u/c_{12}}\,{du}\ll1,
\end{eqnarray*}
By \eqref{555}, for $c_{14}>c_{12}$
$$
\exp\big\{|\,\xi-\eta\,|/c_{14}\tau\big\}\,
{\bf1} \big\{2\sigma\le |\,\xi\,|\le z\big\}\le
\exp\big\{c_{11}\xi^2/c_{12}\sigma^2\big\}\,
{\bf1} \big\{2\sigma\le |\,\xi\,|\le z\big\}
.
$$
Therefore,
$$
{\mathbf E}\,\exp\big\{|\,\xi-\eta\,|/c_{14}\tau\big\}\,
{\bf1} \big\{2\sigma\le |\,\xi\,|\le z\big\}\le
{\mathbf E\exp\big\{v\,W}\big\}\,{\bf1} \big\{ |\,\xi\,|\le z\big\}\ll1
.
$$

By the Cauchy--Bunyakovsky--Schwartz inequality, for $t\in\mathbf{R}$, $z>0$ we have
\begin{equation}
{\mathbf E}\,\exp\big\{t\,|\,\xi\,|+t\,|\,\eta\,|\big\}\,{\bf1} \big\{|\xi|\ge z\big\}\le
\Big({\mathbf E}\,\exp\big\{2t\,|\,\xi\,|\big\}\,{\bf1}\,\big\{|\xi|\ge z\big\}\cdot{\mathbf E}\,\exp\big\{2t\,|\,\eta\,|\big\}\,{\bf1}\,\big\{|\xi|\ge z\big\}\Big)^{1/2},
\label{eq1}
\end{equation}and also
\begin{equation}
{\mathbf E}\,\exp\big\{2t\,|\,\xi\,|\big\}\,{\bf1}\,\big\{|\xi|\ge z\big\}\le
\Big({\mathbf E}\,\exp\big\{4t\,|\,\xi\,|\big\}\cdot{\mathbf P}\,\big\{|\xi|\ge z\big\}\Big)^{1/2}.
\label{eq1tt}
\end{equation}It is clear that
\begin{equation}\label{778}
\exp\big\{4t\,|\,\xi\,|\big\}
\le
\exp\big\{4t\,\xi\big\}+\exp\big\{-4t\,\xi\big\}.
\end{equation}
Applying Lemma~\ref{2.6}, we obtain that for   $z=c_{10}\,\sigma^{2}/\tau$
\begin{equation}\label{0006}
{\mathbf P}\,\big\{|\xi|\ge z\big\}\leqslant 2\, \exp \{-c\sigma^2/\tau^{2}\}.
\end{equation}

Let $0\le |h|\tau\leqslant\frac{1}{2}$.
By~\eqref{6},
\begin{equation}\label{000}
{\mathbf E}\,e^{h\xi}\leqslant  \exp \{h^2\sigma^2\}.
\end{equation}

Applying \eqref{eq1tt}--\eqref{000} with $h=4t=\pm 4/c_{14}\tau$, $z=c_{10}\,\sigma^{2}/\tau$
and choosing a sufficiently large constant $c_{14}$, we obtain
\begin{equation}
{\mathbf E}\,\exp\big\{2t\,|\,\xi\,|\big\}\,{\bf1}\,\big\{|\xi|\ge z\big\}\le
\Big({\mathbf E}\,\big(\exp\big\{4t\,\xi\big\}+\exp\big\{-4t\,\xi\big\}\big)\cdot{\mathbf P}\,\big\{|\xi|\ge z\big\}\Big)^{1/2}\le\sqrt{2}.
\label{eq15}
\end{equation}

It is similarly verified that
\begin{equation}
{\mathbf E}\,\exp\big\{2t\,|\,\eta\,|\big\}\,{\bf1}\,\big\{|\xi|\ge z\big\}\le\sqrt{2}.
\label{eq155}
\end{equation}
Hence,
\begin{equation}
{\mathbf E}\,\exp\big\{|\,\eta\,|/c_{14}\tau+|\,\xi\,|/c_{14}\tau\big\}\,{\bf1} \big\{ |\xi|\ge c_{10}\,\sigma^{2}/\tau\big\}\le2\sqrt{2}.
\label{eq15577}
\end{equation}

Let now $\tau>c_{9}\sigma$. Then

\begin{equation}{\mathbf E}\,\exp\big\{|\,\xi-\eta\,|/c_{14}\tau\big\}
\le
\Big({\mathbf E}\,\exp\big\{2\,|\,\eta\,|/c_{14}\tau\big\}\cdot
{\mathbf E}\,\exp\big\{2\,|\,\xi\,|/c_{14}\tau\big\}\Big)^{1/2}.
\end{equation}
\begin{equation}{\mathbf E}\,\exp\big\{2\,|\,\xi\,|/c_{14}\tau\big\}
\le {\mathbf E}\,\exp\big\{2\,\xi/c_{14}\tau\big\}+
{\mathbf E}\,\exp\big\{-2\,\xi/c_{14}\tau\big\}.
\end{equation}
Applying \eqref{000} with $h=\pm 2/c_{14}\tau$,
and choosing a sufficiently large constant $c_{14}$, we obtain that
\begin{equation}{\mathbf E}\,\exp\big\{|\,\xi-\eta\,|/c_{14}\tau\big\}
\ll1,
\end{equation}completing the proof of the theorem.

Keywords: inequalities, sums of independent random vectors, transport distance estimates, central limit theorem

\end{document}